\documentclass[10pt]{amsart}

\usepackage{amssymb,amsmath}
\usepackage{amsthm,mathabx}
\usepackage{qsymbols}
\usepackage{graphicx}
\usepackage{xcolor}
\usepackage{bbm}

\usepackage[utf8]{inputenc}
\usepackage[T1]{fontenc}
\usepackage{mathtools}
\usepackage{hyperref}

\newtheorem{thm}{Theorem}

\newtheorem{prop}[thm]{Proposition}
\newtheorem{cor}[thm]{Corollary}

\newtheorem*{question}{Question}

\theoremstyle{definition}

\theoremstyle{remark}
\newtheorem{rem}[thm]{Remark}

\newcommand{\R}{\mathbb{R}}

\newcommand{\E}{\mathbb{E}}

\newcommand{\Sd}{\mathbb{S}}

\DeclarePairedDelimiter{\norm}{\lVert}{\rVert} 
\DeclarePairedDelimiter{\abs}{\lvert}{\rvert}
\DeclarePairedDelimiter{\sbra}{[}{]}
\DeclarePairedDelimiter{\pare}{(}{)}

\makeatletter
\let\oldabs\abs
\def\abs{\@ifstar{\oldabs}{\oldabs*}}
\let\oldnorm\norm
\def\norm{\@ifstar{\oldnorm}{\oldnorm*}}
\let\oldsbra\sbra
\def\sbra{\@ifstar{\oldsbra}{\oldsbra*}}
\let\oldpare\pare
\def\pare{\@ifstar{\oldpare}{\oldpare*}}
\makeatother

\begin{document}

\title{Some obstructions to contraction theorems on the half-sphere }

\date{\today}

\author{Max Fathi}
\address{Université Paris Cité and Sorbonne Université, CNRS, Laboratoire Jacques-Louis Lions and Laboratoire de Probabilit\'es, Statistique et Mod\'elisation, F-75013 Paris, France\\
and DMA, École normale supérieure, Université PSL, CNRS, 75005 Paris, France \\
and Institut Universitaire de France (IUF)}
\email{mfathi@lpsm.paris}

\author{Matthieu Fradelizi}
\address{LAMA, Univ Gustave Eiffel, Univ Paris Est Creteil, CNRS, F-77447 Marne-la-Vall\'ee, France.}
\email{matthieu.fradelizi@univ-eiffel.fr}

\author{Nathael Gozlan}
\address{Université Paris Cité, CNRS, MAP5, F-75006 Paris}
\email{nathael.gozlan@u-paris.fr}

\author{Simon Zugmeyer}
\address{UMPA UMR5669, Lyon, France}
\email{simon.zugmeyer@ens-lyon.fr}



\thanks{The second and third named authors are supported by a grant of the Simone and Cino Del Duca foundation. The fourth named author has benefited from a post doctoral position funded by the Simone and Cino Del Duca foundation. This research has been conducted within the FP2M federation (CNRS FR 2036). This research was funded, in whole or in part, by the Agence nationale de la recherche (ANR), Grant ANR-23-CE40-0017 (Project SOCOT) and Grant ANR-23-CE40-0003 (Project CONVIVIALITY). A CC-BY public copyright license has been applied by the authors to the present document and will be applied to all subsequent versions up to the Author Accepted Manuscript arising from this submission, in accordance with the grant’s open access conditions. }
\maketitle

\begin{abstract}Caffarelli's contraction theorem states that probability measures with uniformly log-concave densities on $\R^d$ can be realized as the image of a standard Gaussian measure by a globally Lipschitz transport map. We discuss some counterexamples and obstructions that prevent a similar result from holding on the half-sphere endowed with a uniform measure, answering a question of Beck and Jerison. 
\end{abstract}

\section{Introduction}

The purpose of this work is to discuss counterexamples to certain conjectured global Lipschitz bounds on transport maps on the sphere. The motivation comes from the following theorem of Caffarelli \cite{Caffarelli00}: 

\begin{thm}
Let $\gamma$ be the standard Gaussian measure on $\R^d$, and $\mu$ be a probability measure with density of the form $e^{-V}\gamma$, with $V : \R^n \longrightarrow \R^n \cup\{+\infty\}$ a convex function. Then the Brenier map (or optimal transport map for the quadratic cost) from $\gamma$ onto $\mu$ is $1$-Lipschitz. 
\end{thm}

This theorem has found many applications in probability and analysis, including to correlation inequalities, eigenvalue estimates and rigidity theorems for concentration inequalities. A remarkable feature is that the Lipschitz estimate is \emph{dimension-free}. In particular, it implies that there is a $1$-Lipschitz map from $\gamma$ onto $\mathbf{1}_K\gamma$ whenever $K$ is a convex set. The original proof was based on a maximum principle for the Monge-Amp\`ere PDE that the optimal transport map solves, and more recently there have been proofs based on entropic regularization of optimal transport and properties of the Sinkhorn algorithm \cite{FGP20, CP22}. 

Many comparison theorems in Riemannian geometry for positively curved spaces, such as the Lichnerowicz spectral gap bound or bounds on logarithmic Sobolev constant, and which admit the sphere as optimizer, have counterparts for uniformly log-concave measures, and for which the standard Gaussian is the optimizer (and hence plays the role of the sphere). These Gaussian counterparts (and some rigid versions of them) can often be easily proved using the Caffarelli contraction theorem, just by applying the change of variable (see for example \cite{DPF}). This has led to the question of whether there is a version of the Caffarelli contraction theorem for the sphere, that would imply the aforementioned comparison theorems as consequences. 

Since the optimal transport map between nice continuous probability measures is bijective, one can start by asking for bijective $1$-Lipschitz maps. As pointed out in \cite{Mil18} this extra assumption raises an additional difficulty: positively curved manifolds of the same dimension are not necessarily diffeomorphic. For example, a particular formulation of the above question would be to ask if there always is a bijective $1$-Lipschitz map sending the uniform measure on the $n$-sphere to a given measure on a convex subset of the sphere, with log-concave density. Unfortunately, the answer is trivially no: there cannot be a bijective $1$-Lipschitz map sending the sphere onto a spherical cap, because there must be a point of discontinuity. This has led to raising the same question on the half-sphere instead \cite{BJ21}, and it is that version of the question we shall address here. More generally, \cite{Mil18} also raises the question of existence of $1$-Lipschitz maps between diffeomorphic positively curved manifolds, but here we shall only consider the restricted setting of the half-sphere. 

Let $n\geq 2$ and $\Sd^{n}$ be the unit Euclidean sphere of $\R^{n+1}$, equipped with its geodesic distance $d_{\Sd^n}$. Denote by $\Sd_+^{n} = \Sd^n \cap \{ x \in \R^{n+1} : x_{n+1}\geq 0\}$ the upper hemisphere and by $\sigma$ and $\sigma_+$ the uniform probability measures on  $\Sd^{n}$ and $\Sd_+^{n}$ respectively.

Recall that the $W_2$ distance with respect to the geodesic distance on $\Sd^n$ is defined, for all probability measures $\mu,\nu$ on $\Sd^n$, by
\[
W_2^2(\mu,\nu)= \inf_{\pi \in \Pi(\mu,\nu)} \int d_{\Sd^n}^2(x,y)\,d\pi(x,y)
\]
where $ \Pi(\mu,\nu)$ denotes the set of couplings between $\mu$ and $\nu$. According to McCann's Theorem \cite{McC01}, whenever $\mu$ is absolutely continuous with respect to $\sigma$, there exists a $\mu$ almost surely unique optimal transport map $T$ sending $\mu$ to $\nu$.

The following questions appeared in \cite[Section 4]{BJ21} : \\
Suppose $d\nu = \frac{1}{Z}e^{-V} \mathbf{1}_{K}\,d\sigma_+$ is a probability measure supported on a convex subset $K$ of $\Sd_+^{n}$ with $V: K \to \R$ being a convex function on $K$ :
\begin{itemize}
\item[(a)] Is it true that there exists an injective $1$-Lipschitz map $T:\Sd_+^{n} \to K$ transporting $\sigma_+$ on $\nu$ ?
\item[(b)] Is it true that the $W_2$ optimal transport map between $\sigma_+$ and $\nu$ is $1$-Lipschitz ?
\end{itemize}
The purpose of this note is to show that the answers to these questions are both negative in general. 

In the sequel, a set $K$ is called convex if for any $x,y \in K$, any minimizing geodesic joining $x$ and $y$ entirely lies in $K$. 
A function $f:K \to \R$ is called convex if, for any constant speed minimizing geodesic $\gamma:[0,1] \to K$, it holds 
\[
f(\gamma (t)) \leq (1-t) f(\gamma(0))+tf(\gamma(1)),\qquad \forall t\in [0,1].
\]
By geodesic we will always mean a minimizing geodesic. If $x,y \in \Sd^n$ a geodesic joining $x$ to $y$ is the shorter arc of a great circle passing through $x$ and $y$.

\section{Non-contraction of maps from the half-sphere to uniformly-log-concave measures}
\subsection{Optimal transport maps are not $1$-Lipschitz in general}
First let us see that the question (b) is negative in general. We will denote by $N=(0,\ldots,0,1)$ the north pole of $\Sd^n$.
\begin{prop}\label{prop:W2}
Suppose that $\nu$ is a probability measure on $\Sd^n_+$ with a density of the form $d\nu(x) = f(d(x,N))\,d\sigma_+(x)$, where $f: [0,\pi/2] \to (0,\infty)$ is some positive continuous function. The $W_2$ optimal transport map between $\sigma_+$ and $\nu$ is $1$-Lipschitz if and only if $\nu = \sigma_+$.
\end{prop}

The probability measure $d\nu(x) = \frac{1}{Z}e^{-d^2(x,N)}\mathbf{1}_{\Sd^n_{+}}(x)\,d\sigma_+(x)$ then provides a counterexample to question (b) when $K=\Sd^n_+.$

\begin{rem} Proposition \ref{prop:W2} can be used to rule out assertion (b) for measures supported on convex domains $K$ that are close to the half sphere (see also Section \ref{sec:other target}).
Indeed, let $\nu$ be the probability measure defined above and, for any $0<r<\pi/2$, denote by $\nu_r$ the normalized restriction of $\nu$ to the closed ball of center $N$ and radius $r$. Let $r_k$ be some increasing sequence converging to $\pi/2$. According to \cite[Corollary 5.23]{villani_book} the optimal transport map $T_{r_k}$ from $\sigma_+$ to $\nu_{r_k}$ converges to $T$ (the $\sigma_+$ almost surely unique optimal transport map between $\sigma_+$ and $\nu$) in probability as $k\to \infty$. Extracting a subsequence if necessary, one can assume that pointwise convergence of $T_{r_k}$ to $T$ holds almost everywhere on $\Sd^n_+$. If $T_{r_k}$ was $1$-Lipschitz for all $k$ large enough, then by stability of the class of $1$-Lipschitz maps under pointwise convergence, one would have that $T$ is $1$-Lipschitz too. Therefore there exists a  subsequence  $r_{k'}$ converging to $\pi/2$ along which $T_{r_{k'}}$ is not $1$-Lipschitz. \end{rem}
\proof For notational convenience we deal only with $n=2$ but the reasoning is general.
Let $\psi : \Sd^2 \to \R$ be a $c$-concave maximizer of the dual problem for the quadratic transport between $\sigma_+$ and $\nu$ :
\[
\frac{1}{2}W_2^2(\sigma_+,\nu) = \int \psi\,d\sigma_+ + \int \psi^c\,d\nu,
\]
where, for all $x,y\in \Sd^2$, $c(x,y)=\frac{d^2(x,y)}{2}$, $\psi^c(y)=\inf_x\left\{c(x,y) - \psi (x)\right\}$ and $c$-concavity of $\psi$ means that $\psi^{cc}=\psi$. 
For any $\theta \in \R$, denote by $R_{\theta}$ the rotation matrix 
\[
\left[\begin{array}{ccc}  \cos \theta& -\sin \theta  & 0   \\ \sin \theta  & \cos \theta  & 0  \\  0& 0   & 1 \end{array}\right].
\] 
Since $(R_\theta)_\#\sigma_+ = \sigma_+$ and $(R_\theta)_\#\nu=\nu$ for all $\theta \in \R$, one sees that if $\psi$ is a dual maximizer, then $\psi\circ R_\theta$ also (note that $(\psi\circ R_\theta)^c = \psi^c\circ R_{-\theta}$). Moreover, if $\Theta$ is a random variable uniformly distributed on $[-\pi,\pi]$, then denoting by $\bar{\psi}(x) = \E[\psi (R_{\Theta}x)]$, $x \in \Sd^2$, one sees that $\int \bar{\psi}\,d\sigma_+ = \int \psi\,d\sigma_+$ and 
\[
\bar{\psi}^c(y) = \inf_{x\in \Sd^2} \E[c(x,y) - \psi(R_\Theta x)] \geq \E[\psi^c(R_{-\Theta}y)].
\]
Thus, 
\[
\int \bar{\psi}^c\,d\nu \geq \int \psi^c\,d\nu.
\]
Therefore, $\bar{\psi}$ is a dual maximizer and is rotationally invariant : $\bar{\psi}\circ R_\theta = \bar{\psi}$ for any $\theta$. Setting $\varphi = \bar{\psi}^{cc}$, it is easy to check that $\varphi$ is rotationally invariant, $c$-concave and a dual maximizer. According to McCann's theorem, the optimal transport map between $\sigma_+$ and $\nu$ is given, for $\sigma_+$ almost all $x \in \Sd_+^2$, by
\[
T(x) = \exp_x(- \nabla \varphi(x)).
\]
This optimal transport map $T$ is almost surely unique (and does not depend on the particular dual potential).
 Moreover, since $\varphi$ is invariant under all $R_\theta's$, there exists $h:[0,\pi] \to \R$ such that $\varphi(x)=h(d(x,N))$. Thus for almost every $x\in \Sd_+^2$, $T(x)$ belongs to the geodesic connecting $x$ to $N$. 

Now let us assume by contradiction that $T:\Sd_+^2 \to\Sd_+^2$ is $1$-Lipschitz.
By rotational invariance of $\nu$ and $\sigma_+$, uniqueness and continuity of $T$, one easily gets that the equality
\[
R_{\theta}\circ T \circ R_{-\theta} = T
\]
holds everywhere on $\Sd^2_+$.
So $R_{\theta}(T(N))= T(N)$ for all $\theta$ and so $T(N)=N$.
Fix some $x$ such that $d(x,N)=\pi/2$ and for all $t\in [0,\pi/2]$ let $x_t$ be the unique point on the geodesic joining $N$ to $x$ such that $d(N,x_t)=t$. Then, there exists a function $r:[0,\pi/2] \to [0,\pi/2]$ such that
\[
T(x_t) = x_{r(t)},\qquad \forall t\in [0,\pi/2].
\]
Moreover,
\[
|r(t) - r(s)| = d(x_{r(t)},x_{r(s)}) = d(T(x_t),T(x_s)) \leq d(x_t,x_s)=|t-s|
\]
and so $r$ is $1$-Lipschitz. By rotational invariance, the function $r$ is in fact independent of $x$. Therefore, 
\[
T(\Sd_+^2) = \{ y \in \Sd_+^2 : d(N,y) \in r([0,\pi/2])\}.
\]
Since $T(\Sd_+^2)=\Sd_+^2$, $r$ is necessarily surjective. Therefore, $r: [0,\pi/2] \to [0, \pi/2]$ is a $1$-Lipschitz and surjective map. The only possibility is that $r (x) = x$ for all $x \in [0,\pi/2]$ or $r(x) = \pi/2 - x$ for all $x \in [0,\pi/2]$. Indeed, if $a,b \in [0,\pi/2]$ are such that $r(a)=0$ and $r(b)=\pi/2$, then $r$ being $1$-Lipschitz, it holds $\pi/2\leq |a-b|$ and so $|a-b|=\pi/2$. Suppose that $a=0$ and $b=\pi/2$. Since $r$ is $1$-Lipschitz, the function $k(x)=x\mapsto r(x)-x$ is non-increasing and such that $k(0)=k(1)=0$. Therefore $k(x)=0$ for all $x \in [0,\pi/2]$. 
If $a=\pi/2$ and $b=0$, one sees similarly that $r(x)= \pi/2-x$ for all $x \in [0,\pi/2].$ Since we know that $r(0)=0$ the only possibility is $r=Id$ and so $T = Id$.
\endproof

\subsection{Bijective transport maps are not $1$-Lipschitz in general}

We now want to show that the answer to question (a) is also negative. The idea is that if the target measure has the half-sphere as support, then such a map would preserve the boundary, and the $1$-Lipschitz constraint would then force its value along great circles joining opposite points on the boundary. The first part of this reasoning is given by the following result:
\begin{prop}\label{prop:isom1}
Assume that $T: \Sd^n_+ \to \Sd^n_+$  is a surjective $1$-Lipschitz map. Denote by $C = \Sd^n_+ \cap\{x_{n+1}=0\}$ and assume that $T(C)=C$. Then $T$ is an isometry fixing $N$.
\end{prop}

The following result is due to Valentine \cite[Lemma 3]{Val45}.
\begin{thm} \label{thm:Valentine}
Let $T : A \to \Sd^{n}$ be defined on a compact subset $A$ of $\Sd^{n}$ with Lipschitz constant $L \leq 1$ ; if the image of $T$ is not contained in a closed hemisphere, then $T$ is an isometry. 
\end{thm}
A proof can be found in e.g Appendix D of \cite{Brudnyi12} or in \S 13 of \cite{WW75}.

\proof[Proof of Proposition \ref{prop:isom1}]
Let $N=(0,\ldots, 0,1) \in \Sd^n$ be the north pole. For any $x \in C$, it holds $d(T(N), T(x)) \leq d(N,x)=\pi/2$. Since $T(C)=C$, one gets that $d(T(N), x') \leq \pi/2$ for all $x' \in C$. This forces $N=T(N)$.
If $(x_t)_{t\in [0,\pi/2]}$ is a geodesic such that $x_0 \in C$ and $x_{\pi/2}=N$, then $y_t = T(x_t)$ satisfies $d(y_s,y_t) \leq d(x_s,x_t)=|t-s|$ and $d(y_0,y_{\pi/2}) = \pi/2$. Therefore, $(y_t)_{t\in [0,\pi/2]}$ is the geodesic connecting $T(x_0)$ to $N$. Therefore the map $T$ is completely determined by its restriction to $C.$ Since $C = \Sd^{n-1} \times\{0\}$ and $T(C)=C$, one sees that $T$ induces a $1$-Lipschitz surjective map from $\Sd^{n-1}$ to itself. 
According to Valentine's Theorem \ref{thm:Valentine}, it is an isometry of $\Sd^{n-1}$ and so is the restriction of an isometry of $\R^n$ to $\Sd^{n-1}$. 
In other words, there is an $n\times n$ orthogonal matrix $A$ such that
\[
T((z,0)) = (Az,0),\qquad \forall z\in \Sd^{n-1}.
\]
Now, if $x = (z,x_{n+1}) \in \Sd^n$, then $d(x,N) = \arccos(x_{n+1})$ and $x$ is on the geodesic connecting $N$ to $(z/|z|,0)$. So, $T(x)$ is the unique point at distance $ \arccos(x_{n+1})$ from $N$ on the geodesic going from $N$ to $(A(z/|z|),0)$. Since $|Az|=|z|$, one gets $T(x) = (Az,x_{n+1})$.
\endproof

\begin{cor}\label{Cor:1}
If $T: \Sd^n_+ \to \Sd^n_+$  is a bijective $1$-Lipschitz map, then $T$ is an isometry fixing $N$.
In particular the image of $\sigma_+$ under $T$ is $\sigma_+$.
\end{cor}
In particular, if $\sigma_+ \neq \nu$ and $K=\Sd^n_{+}$ then there is no injective and $1$-Lipschitz map transforming $\sigma_+$ into $\nu$.
\proof
Denote by $\Sd^n_{++} = \{x \in \Sd^n : x_{n+1}>0\}$. Since $T:\Sd^n_+ \to \Sd^n_+$ is an homeomorphism, it is clear that $T(\Sd^n_{++})$ is an open subset of $\Sd^n_+$ (for the relative topology of $\Sd^n_+$). Let us show that $T(\Sd^n_{++})$ is also open in $\Sd^n$.
Let $B= \{x \in \R^n : |x| <1\}$ and consider the map $S : \Sd^n_+ \to \bar{B} : x \mapsto (x_1,\ldots,x_n)$ and the map $U : \bar{B} \to \Sd^n_+ : y \mapsto (y, \sqrt{1- |y|^2})$. The function $f:\bar{B} \to \bar{B}$ defined by $f = S\circ T \circ U$ is continuous and bijective. By Brouwer's invariance of domain theorem, $f(B)$ is an open subset of $\R^n$.  Therefore, $T(\Sd^n_{++}) = U\circ f (B)$ is an open subset of $\Sd^n$. Denoting as before $C = \Sd^n_+ \cap\{x_{n+1}=0\}$, one thus gets that $T(\Sd^n_{++}) \cap C = \emptyset$. Since $\Sd^n_+ = T(C) \cup T(\Sd^n_{++})$, one concludes that $C \subset T(C)$. Denoting by $\mathcal{H}_{n-1}$ the $n-1$ dimensional Hausdorff measure on $\Sd^n$, one gets 
\[
\mathcal{H}_{n-1}(T(C)) \leq \mathcal{H}_{n-1}(C),
\]
since $T$ is $1$-Lipschitz. Thus $\mathcal{H}_{n-1}(T(C) \setminus C)=0 $. In particular, $C$ is dense in $T(C)$. Since $C$ is closed, one gets $T(C)=C.$ Using Proposition \ref{prop:isom1} completes the proof.
\endproof

A consequence is the following strong improvement of Proposition \ref{prop:W2}:
\begin{cor}\label{cor:W21Lip}
Suppose that $\nu$ has a density of the form $d\nu(x) = g(x)\,d\sigma_+(x)$, where $g: \Sd^n_+ \to (0,\infty)$ is some measurable function such that $a\leq g\leq b$ for some $0<a<b$.  Then the $W_2$ optimal transport map $T$ between $\sigma_+$ and $\nu$ is $1$-Lipschitz if and only if $\nu = \sigma_+$.
\end{cor}
\proof
Let $S$ be the optimal transport map sending $\nu$ onto $\sigma^+$ given by Mc Cann's theorem and let us justify that $S$ is continuous.  Let $\bar{\nu} = \frac{1}{2} \nu + \frac{1}{2}\nu_-$, where $\nu_-$ is the image of $\nu$ under the map 
\[
u :\Sd^n \to \Sd^n : x \mapsto (x_1,\ldots,x_{n-1},-x_n)
\] 
and denote by $\bar{S}$ the optimal transport map sending $\bar{\nu}$ on $\sigma$. 
Let us admit for a moment that 
\begin{equation}\label{eq:barS}
\bar{S}(x) = S(x),\qquad \forall x \in \Sd^n_+.
\end{equation}
According to \cite[Theorem 2.4]{Loeper11}, the transport map $\bar{S}$ is continuous, and so \eqref{eq:barS} yields that $S$ is continuous on $\Sd^n_+$.  Let $T$ be the $W_2$ optimal map transporting $\sigma_+$ onto $\nu$  and assume that $T$ is $1$-Lipschitz. Since $\nu$ has full support, $T$ is surjective.  According to \cite[Corollary 10]{McC01}, $S$ and $T$ satisfy $T\circ S(x)=x$ and $S\circ T(x)=x$ for $\sigma_+$ almost every $x \in \Sd^n_+$. Since $S$ is continuous, these identities actually holds for all $x \in \Sd^n_+$. Thus $T$ is injective. Applying Corollary \ref{Cor:1}, one gets that $T$ is an isometry fixing $N$ and so $\nu=\sigma_+$.
Now let us show \eqref{eq:barS}. Since $\bar{\nu}$ and $\sigma$ are invariant under the isometry $u$, by uniqueness of the optimal transport map, $T$ satisfies $T = u \circ T \circ u$, $\sigma$ a.e and, by continuity of $T$, everywhere.
Let us show that for all $x \in \Sd^n_{++}$, $T(x) \in \Sd^n_{+}$, which will imply \eqref{eq:barS} by continuity of $T$.
Suppose, by contradiction, that some $x \in \Sd^n_{++}$ is such that $T(x)  \in \Sd^n_{--} := \Sd^n \setminus \Sd^n_+.$ Then $x'= u(x) \in \Sd^n_{--}$ is such that $T(x') = u(T(x)) \in \Sd^n_{++}$. 
By cyclical monotonicity, it holds 
\begin{equation}\label{eq:cyclical}
d^2(x,T(x)) + d^2(x',T(x')) \leq d^2(x,T(x'))+d^2(x',T(x)).
\end{equation}
But 
\[
x\cdot T(x') = x\cdot  u(T(x)) = \sum_{i=1}^{n-1} x_i(T(x))_i -x_n(T(x))_n >\sum_{i=1}^{n-1} x_i(T(x))_i +x_n(T(x))_n= x\cdot T(x).
\]
The function $\arccos$ being decreasing, one gets that
\[
d^2(x,T(x')) <d^2(x,T(x)) 
\]
and similarly $d^2(x',T(x)) <d^2(x',T(x'))$, contradicting \eqref{eq:cyclical}. This completes the proof.
\endproof

\begin{rem}Let $\nu$ be as Corollary \ref{cor:W21Lip} and let, as in the preceding proof,  $\bar{S}$ be the optimal transport map sending $\bar{\nu}$ onto $\sigma$ and $\bar{T}$ be the transport map sending $\sigma$ on $\bar{\nu}$. Reasoning as above, we see that $\bar{S}$ and $\bar{T}$ are both continuous (using again \cite[Theorem 2.4]{Loeper11}), are inverse of each other, and such that $T(\Sd^n_+)=\Sd^n_+$ and $T(\Sd^n_-)=\Sd^n_-$. Thus 
\[
T(C) = T(\Sd^n_+ \cap \Sd^n_-) = T(\Sd^n_+) \cap T(\Sd^n_-)=\Sd^n_+ \cap \Sd^n_-=C.
\]
So, to prove Corollary \ref{cor:W21Lip}, one could have directly use Proposition \ref{prop:isom1} instead of Corollary \ref{Cor:1}.
\end{rem}

\begin{rem}\label{rem_unbounded_lip} Let us give some more quantitative version of Corollary \ref{cor:W21Lip}.
For all $\varepsilon>0$, consider a probability measure of the form
\[
d\nu_\varepsilon(x)= \frac{1}{Z_\varepsilon}e^{-\frac{1}{\varepsilon} V(x)}\,d\sigma_+(x)
\]
with $V$ some bounded continuous function on $\Sd^n_{++}$ admitting $N$ as unique minimizer. 
As $\varepsilon \to 0$, $\nu_\varepsilon \to \delta_N$ for the weak topology.
Let $m_\varepsilon = \int d(x,N)\,d\nu_\varepsilon(x)$  and observe that  $m_\varepsilon \to 0$ as $\varepsilon \to 0$. For  $\varepsilon$ small enough, let $r_\varepsilon>0$ be such that
\[
\nu_\varepsilon(B(N,r_\varepsilon))=1-\sqrt{m_\varepsilon},
\]
where $B(N,r_\varepsilon)$ is the closed ball of radius $r_\varepsilon$ centered at $N$.
By Markov's inequality, $m_\varepsilon \geq r_{\varepsilon} \sqrt{m_\varepsilon}$, and $m_\varepsilon \to 0$ as $\varepsilon \to 0$. Therefore $r_\varepsilon \to 0$ as $\varepsilon \to 0$. Let $T_\varepsilon$ be the optimal transport map sending $\sigma_+$ to $\nu_\varepsilon$ and define $K_\varepsilon = T_\varepsilon^{-1}(B(N,r_\varepsilon))$. As explained in the preceding remark, $T_\varepsilon$ is a continuous map such that $T_\varepsilon(C)=C$. Thus $K_\varepsilon$ is a closed subset of $\Sd^n_+$ such that $K_\varepsilon \cap C = \emptyset$ and $\sigma_+(K_\varepsilon)=1-\sqrt{m_\varepsilon}$. Denote by $R_\varepsilon= \sup_{x \in K_\varepsilon} d(x,N)$, which satisfies $R_\varepsilon<\pi/2.$ Since $\sigma_+(B(N,R_\varepsilon) \to 1$ as $\varepsilon \to 0$, one concludes that $R_\varepsilon\to \pi/2$.
If $x_\varepsilon$ is some point in $K_\varepsilon$ such that $d(x_\varepsilon,N)=R_\varepsilon$ and $x'_\varepsilon$ is the unique point on $C$ such that $x_\varepsilon$ belongs to the geodesic joining $x'_\varepsilon$ to $N$, then one gets
\[
d(T_\varepsilon(x_\varepsilon), T_\varepsilon(x'_\varepsilon)) \geq d(B(N,r_\varepsilon),C)= \pi/2 - r_\varepsilon
\]
and so
\[
\frac{d(T_\varepsilon(x_\varepsilon), T_\varepsilon(x'_\varepsilon))}{d(x_\varepsilon,x'_\varepsilon)} \geq \frac{\pi/2 - r_\varepsilon}{\pi/2 - R_\varepsilon}
\]
and so 
\[
\mathrm{Lip}(T_\varepsilon) \geq \frac{\pi/2 - r_\varepsilon}{\pi/2 - R_\varepsilon} \to +\infty
\]
as $\varepsilon \to 0.$
\end{rem}

Finally, let us recall the following well known fact that the classes of $1$-Lipschitz functions from $\Sd^n$ to $\Sd^n$ with respect to geodesic or Euclidean distances on the sphere coincide. The conclusion of Corollary \ref{cor:W21Lip} is thus also true if $\Sd^n$ is equipped with the ambient Euclidean distance of $\R^{n+1}$.

\begin{prop}
A map $T : A \to \Sd^n$ defined on a subset $A$ of $\Sd^n$ is $1$-Lipschitz with respect to the geodesic distance if and only if it is $1$-Lipschitz with respect to the Euclidean norm $\|\,\cdot\,\|$ on $\R^{n+1}$.
\end{prop}
\proof
Since $d_{\Sd^n}(x,y) = \arccos (x\cdot y)$ and $\cos$ is decreasing on $[0,\pi]$ one has
\[
d_{\Sd^n}(T(x),T(y)) \leq d_{\Sd^n}(x,y)  \quad \Rightarrow  \quad T(x)\cdot T(y) \geq x\cdot y.
\]
Since $x,T(x), y,T(y)$ belong to $\Sd^n$, it holds $\| T(x) - T(y)\|^2 = 2(1-T(x) \cdot T(y))$ and $\| x -y\|^2 = 2(1-x \cdot y)$, which completes the proof.
\endproof

\begin{rem}The conclusion is false for Lipschitz constants $<1$. 
Namely, let $n=1$ and $A = \Sd^1 \cap \R^+\times \R$. For $\theta \in [-\pi/2,\pi/2]$ define $u(\theta) = (\cos(\theta),\sin(\theta))$ and consider the map $T: A \to \Sd^1$ defined by $T(u(\theta)) = u(\theta/2)$, for all $\theta \in [-\pi/2,\pi/2]$. The map $T$ is $1/2$-Lipschitz with respect to the geodesic distance. Indeed, $d_{\Sd^1}(u(\theta_1),u(\theta_2)) = |\theta_2-\theta_1|$, for all $\theta_1,\theta_2 \in [-\pi/2,\pi/2].$ Observe however that $T$ is not $1/2$-Lipschitz for the Euclidean norm on $\R^2$. Indeed, $\|T(\pi/2) - T(-\pi/2)\| = \|u(\pi/4)- u(-\pi/4)\|= \sqrt{2}$ whereas, $\|u(\pi/2)-u(-\pi/2)\| = 2$, so that $\|T(\pi/2) - T(-\pi/2)\| =\sqrt{2} > 1 = \frac{1}{2}\|u(\pi/2)-u(-\pi/2)\|$.
\end{rem}

\section{Other target domains}\label{sec:other target}

Corollaries \ref{Cor:1} and \ref{cor:W21Lip} rule out existence of $1$-Lipschitz transport maps when the domain of the target measure is the full half-sphere. Of course, when the target domain is smaller, it is possible to have a $1$-Lipschitz transport map. However, by mimicking the argument made above, we can see that if a Lipschitz transport map exists, then the target measure cannot be too concentrated away from its boundary, in the following sense: 

\begin{prop}
Let $\nu$ be a probability measure supported on a closed subset $D \subset \Sd_+^{n}$, and assume that there is a bijective $L$-Lipschitz transport map $T$ sending $\sigma_+$ onto $\mu$. Then 
\[
\nu\left( \Sd^n_+ \setminus (\partial D)^{r}\right)\leq 2\exp\left(-(n-1)r^2/(2L^2)\right),\qquad \forall r < \frac{\pi L}{2}.
\]
\end{prop}
Here $A^r$ denotes the open $r$-neighborhood of the set $A$, with respect to the distance on the sphere, that is
$$A^r := \{y \in \Sd_+^{n}; d(y, A) < r\}.$$ 
Therefore the mass away from the boundary $\partial D$ cannot be too large. Note that this easily implies the lack of uniform estimate on Lipschitz norms of transport maps onto log-concave densities on the half-sphere, as discussed in Remark \ref{rem_unbounded_lip}.

\begin{proof}
Since $T$ is bijective and continuous on a compact set, it is a homeomorphism, and as a consequence of the invariance of domain theorem, it maps $C = \Sd^n_+ \cap\{x_{n+1}=0\}$ onto $\partial D$. As it is $L$-Lipschitz, any point from $B(N,R)^c$ (the complement of the closed ball of radius $R$ around the north pole, with respect to the spherical distance on $\Sd^n_+$ is mapped to a point at a distance at most $L(\pi/2 - R)$ from $\partial D$. 
Thus, taking $R = \pi/2 - r/L$, with $r< L \pi/2$, it holds 
\[
T(B(N,\pi/2 - r/L)^c) \subset (\partial D)^{r}
\]
and so
\[
T^{-1} \left( D \setminus (\partial D)^{r}\right) \subset B(N,\pi/2 - r/L).
\]
Therefore 
\[
\nu\left( D \setminus (\partial D)^{r} \right) = \sigma_+\left( T^{-1} \left( D \setminus (\partial D)^{r}\right) \right) \leq \sigma_+(B(N, \pi/2 - r/L)) .
\]
Finally, recall that the uniform measure $\sigma$ on $\Sd^n$ satisfies the following Gaussian concentration inequality (see e.g \cite{Ledoux01}): for all $A \subset \Sd^n$ such that $\sigma (A) \geq 1/2$, 
\[
\sigma(\Sd^n \setminus A^t) \leq e^{-(n-1)t^2/2}, \qquad \forall t\geq 0.
\]
Applying this inequality to $A = \Sd^n_-$ and $t = r/L$ gives that
\[
\sigma_+(B(N, \pi/2 - r/L)) = 2 \sigma (\Sd^n \setminus \left(\Sd^n_-\right)^{r/L}) \leq 2\exp\left(-(n-1)r^2/(2L^2)\right)
\]
which completes the proof.
\end{proof}

Finally, we would like to stress that the following sub-case of the question of Beck and Jerison does not seem to be ruled out by our results. 

\begin{question}
Is there a $1$-Lipschitz transport map from $\sigma_+$ onto $\frac{1}{\sigma_+(K)}\mathbf{1}_K\sigma_+$ for any convex subset $K$ of $\Sd^n_+$ with positive volume ? 
\end{question}

\bibliographystyle{alpha}

\end{document}